\documentclass{amsart}

\usepackage{amssymb}
\usepackage{xcolor}
\usepackage{bbm}

\newtheorem{thm}{Theorem}[section]
\theoremstyle{definition}

\newtheorem{prop}[thm]{Proposition}
\newtheorem{coro}[thm]{Corollary}
\newtheorem{defn}[thm]{Definition}
\newtheorem{lem}[thm]{Lemma}

\theoremstyle{remark}

\allowdisplaybreaks

\begin{document}

\title[Weighted Bergman projection in radial trees]{Weighted Bergman projection in radial trees}


\author{Nikolaos Chalmoukis}
\address{Department of Mathematics and Applications, University of Milano-Bicocca, Piazza dell'Ateneo Nuovo, 20126 Milano, Italy}
\email{nikolaos.chalmoukis@unimib.it}

\author{Matteo Monti}
\address{Department of Mathematics and Applications, University of Milano-Bicocca, Piazza dell'Ateneo Nuovo, 20126 Milano, Italy}
\email{matteo.monti@unimib.it}


\begin{abstract}
We study the Bergman projection $P_\omega$ on a radial tree with respect to weighted harmonic Bergman spaces $A^p_\omega$, where $\omega$ is a positive radial weight. We prove that for all such weights and $1<p<\infty$, the Bergman projection $P_\omega$ is bounded from $L^p(\omega)$ to $A^p_\omega$ and $L^\infty \to \mathcal{B}$, where $\mathcal{B}$ is the Bloch space of the tree. We also provide a characterization of the weights $\omega$ for which the positive Bergman projection $P_\omega^+$ is bounded from $L^p(\omega)$ to $L^p(\omega)$. 
Finally, in the appendix, we provide an elementary proof of the characterization of the weights $\Omega \in L^1(\mathbb{D})$ for which the Bergman projection from $L^p(\Omega)$ to the holomorphic Bergman space $A^p_\Omega(\mathbb{D})$.
\end{abstract}


\maketitle

\section{Introduction}

In the recent years, there has been a growing interest in the study of problems related to function spaces and harmonic analysis on combinatorial trees. At times the interest stemmed from an attempt to formulate toy versions of problems in hyperbolic spaces, with the intend to transfer the results back to the continuum. A classical application of this principle in the domain of function spaces is the article \cite{Arcozzi2008}, while some more recent examples are \cite{Chalmoukis2026} and \cite{McCourt2026}. 
It is also the case that one would like to develop a theory of function spaces on trees parallel to the one on hyperbolic spaces in order to observe where the two theories converge and where they diverge. Such an exercise offers a deeper understanding of the phenomena that govern both worlds and of the reasons that the behaviour of certain objects is different. 
The scope of this article falls in the second category. We will try to develop a part of the theory of weighted Bergman spaces on radially symmetric trees that is missing from the literature and compare it with the theory in the unit disc.

 We will denote by $X $ a rooted tree, that is a locally finite connected graph with no cycles. The root vertex is a distinguished vertex of $X $ denoted by $o $. With an abuse of notation we will identify the graph $X $ with the set of its vertices. 
 We write $u \sim v $  for $ u,v \in X $ if there is an edge connecting $u $ and $v$. In such a case, we say that $u $ and $v $ are neighbours. 
 The tree $X $ is naturally endowed with the edge counting distance $d $. The distance of any vertex $ v $ from the origin is called the level of $v $ and denoted by $ |v| $. Formally, $ |v|:= d(v,o) $. We call the function defined by $q(v) := \# \{ u\in X : u\sim v \} - 1 $ for $v\neq o $ , $q(o):= \# \{ u\in X : o \sim u \} $  the branching function of $X$, where $ \# $ denotes the cardinality of a set.
 Throughout this article we assume that
 \begin{equation}\tag{$A_0$} \label{hyp:H0}
	 2 \leq q(v) < + \infty,  
 \end{equation}
 for every $ v\in X $, and that $q $ is radial, that is there exists $\overline{q}: \mathbb{N}_0\to \mathbb{N}_0 := \mathbb{N} \cup \{ 0 \} $ such that
 \begin{equation} \tag{$A_1 $ }
	 q(v) = \overline{q}(|v|)
 \end{equation}
 for $ v \in X$. With a small abuse of notation we shall use $q$ in place of $\overline q$.  
 Notice that we do not require $q $ to be a bounded function on the tree. Given $ \alpha > 1 $ we can define a natural finite measure on the tree as follows 
 \[ \mu_\alpha(x) = \Big( \prod_{\ell=0}^{|x|-1} \frac{1}{q(\ell)} \Big)^\alpha. \]
Here and in the following we adopt the convention that an empty product is equal to $1$ and an empty sum equals $0 $.

We will now introduce the function spaces that we are mainly interested in on a radially symmetric tree. A function $ f: X \to \mathbb{C} $ is called harmonic if 
\[ f(v) = \begin{cases} {\displaystyle \frac{1}{q(v)+1}\sum_{u \sim v} f(u)}, & v\neq o  \\
{\displaystyle \frac{1}{q(o)} \sum_{u \sim o } f(u)}, & v=o   \end{cases}.   \]
Let $1\leq p<\infty $ and $ \omega \in L^1(X) $ a positive function. The space 
 \[ A_\omega^p = \{ f\in L^p(\omega) : f \, \,\, \text{is harmonic on}\,\, X \} \]
is called harmonic Bergman space of exponent $p$ with respect to $\omega. $ 
It is immediate to see that $ A_\omega^p $ is a Banach space and in particular, for $p=2 $, it is a closed subspace of the Hilbert space $ L^2(\omega) $, so that there exists an orthogonal projection $ P_\omega : L^2(\omega) \to A^2_\omega. $  From now on we assume $\omega $ to be a radial weight, i.e. there exists a function $\overline{\omega}: \mathbb{N}_0\to (0,+\infty) $ such that 
\begin{equation}
	\tag{$A_3$}
		\omega(v) = \overline{\omega}(|v|), \quad v\in X.
\end{equation}

The weighted Bergman spaces were first introduced in \cite{Picardello2016} and the study of the boundedness of the operator $P_\omega : L^p(\omega) \to A^p_\omega $ was initiated in \cite{DeMari2024}. In that work the authors prove that  for homogeneous trees, which correspond to a constant function $q $, the operator  $P_{\mu_\alpha} : L^p(\mu_\alpha) \to A^p_{\mu_\alpha}  $ is bounded for all $ 1 < p < \infty $ .  
The result was extended to non homogeneous radial trees in \cite{CondeAlonso2026} for the same class of measures.

These developments are a reflection of the study of weighted Bergman spaces in the unit disc. 
In fact, if we denote by $ \mathbb{D} = \{ z\in \mathbb{C} : |z| < 1 \} $ the unit disc, and by $ dA(z) = \frac{dxdy}{\pi} $ the normalized area measure, we can define for a positive continuous function $ \Omega \in L^1(\mathbb{D}, dA) $ and $ 1 \leq p < +\infty $  the weighted Bergman spaces $ A_\Omega^p(\mathbb{D}) $ as the space of holomorphic functions $ f: \mathbb{D} \to \mathbb{C} $ such that  
\[ \int_\mathbb{D} |f(z)|^p \Omega(z)dA(z) < + \infty. \]

As in the case of trees, the space $ A_\Omega^p(\mathbb{D}) $ is a closed subspace of the Hilbert space $ L^2(\Omega dA) $ and hence there exists an orthogonal projection $ \Pi_\Omega : L^2(\Omega dA) \to A_\Omega^2(\mathbb{D}) $. For the case of standard weights of the form $ \Omega_\alpha(z) = (1-|z|^2)^\alpha $ with $ \alpha > -1 $, it is known that the operator $ \Pi_{\Omega_\alpha} : L^p(\Omega_\alpha dA) \to A^p_{\Omega_\alpha}(\mathbb{D}) $ is bounded for all $ 1 < p < \infty $ as it was proved in \cite{Forelli1974}. 

The more general class of radial symmetric weights is already quite interesting. A weight $\Omega $ is radially symmetric if $ \Omega(z) = \Omega(|z|)$ for $z\in \mathbb{D}. $ 
Dostani\'c \cite{Dostanic2004, Dostanic2009} observed that a necessary condition for the boundedness $ \Pi_{\Omega} : L^p(\Omega dA) \to A^p_\Omega $ for a fixed $p$, $1<p<\infty $, is the so-called Dostani\'c condition, i.e. there exists $C > 0 $ such that  
\begin{equation}\label{cond:Dostanic}
	\tag{$D_p $} 
	 \Big( \int_0^1 \Omega(t) t^{np+1} dt  \Big)^{\frac{1}{p}} \Big( \int_0^1 \Omega(t) t ^{np'+1} \Big)^{\frac{1}{p'}} \leq C \int_0^1 \Omega(t) t ^{2n+1} dt, \,\, \text{for all}\,\, n\in \mathbb{N}
\end{equation}
where $p' $ is the H\"older conjugate exponent of $p. $ 

Furthermore the exponentially decreasing weights 
\[ \Omega(z) = \exp \bigg( - \frac{1}{(1-|z|^2)^\alpha} \bigg), \qquad |z|<1, \alpha > 0 \]
do not satisfy condition \eqref{cond:Dostanic} for any $ p \neq 2 $ .

A systematic study of the boundedness properties of the Bergman projection for radial weights was undertaken by Pel\'aez and R\"atty\"a \cite{Pelaez2021}. 
They consider the class of upper doubling weights $ \widehat{\mathcal{D}} $. A radial weight $\Omega $ belongs to $\widehat{\mathcal{D}} $ if there exists $C>1$ such that 
\begin{equation}
	\label{cond:upperD} \tag{UD}
	\int_r^1 \Omega(t) dt \leq C \int_{\frac{1+r}{2}}^1\Omega(t) dt, \,\, \text{for all}\,\, r\in (0,1).
\end{equation}
By \cite[Theorem 7]{Pelaez2021} it follows that if $\Omega \in \widehat{\mathcal{D}} $ then $ \Pi_\Omega : L^p(\Omega dA) \to A^p_\Omega $ is bounded for all $ p \in(1,\infty). $  
Recently, Li and Wang \cite{Li2026} closed the loop by showing that if $ \Omega $ satisfies the Dostani\'c condition for any $ p \in (1,\infty)\setminus \{ 2 \} $ then $\Omega \in \widehat{\mathcal{D}}. $  Therefore proving the following theorem 
\begin{thm}
	Let $ \Omega $ a positive radial weight in $ L^1(dA) $ then the following are equivalent.
	\begin{itemize}
		\item[(a)] $ \Pi_\Omega : L^p(\Omega dA) \to A^p_\Omega $ is bounded for some $ p\in (1,\infty) \setminus \{ 2 \} $.
		\item[(b)] $ \Pi_\Omega : L^p(\Omega dA) \to A^p_\Omega $ is bounded for all $ p \in (1,\infty) $ .
		\item[(c)] $\Omega \in \widehat{\mathcal{D}} $ .
	\end{itemize}
\end{thm}
Their proof in \cite{Li2026} quite intricate, therefore we deemed useful to include a short elementary proof in the Appendix (see Section \ref{sec: Appendix}).

Furthermore, if one considers the integral operator $ \Pi_\Omega^+ $ which has as kernel the absolute value of the kernel of the Bergman projection $ \Pi_\Omega, $ Pel\'aez and R\"atty\"a \cite[Theorem 9]{Pelaez2021} observed that it has different boundedness properties with respect to the Bergman projection, contrary to the situation for the standard weights $ \Omega_\alpha $ where it is well known that both operators are bounded for all $p \in (1,\infty) $ .
More precisely, the authors introduce the class of lower doubling measures $ \check{\mathcal{D}} $ which comprises radial weights $\Omega $ such that there exists $ K>1 $ and $ C> 1 $ such that 
\begin{equation}
	\label{cond:LD} \tag{LD} 
	\int_r^1 \Omega(t) dt \geq C \int _{1-\frac{1-r}{K}}^1 \Omega(t) dt, \,\, \text{for all}\,\, r \in (0,1).
\end{equation}
Finally, we define the class $ \mathcal{D} = \widehat{\mathcal{D}}\cap \check{\mathcal{D}} $ . Then \cite[Theorem 9]{Pelaez2021} says that  if $\Omega \in \widehat{\mathcal{D}} $ and $ p\in (1,\infty) $ the positive Bergman projection $\Pi_\Omega^+ $ is bounded on $ L^p(\Omega dA) $ if and only if $ \Omega \in \mathcal{D}. $ Combined with the observation that the boundedness of $\Pi_\Omega^+ $ implies the boundedness of $\Pi_\Omega $ we obtain that $ \Pi_\Omega^+ : L^p(\Omega dA) \to L^p(\Omega dA) $ is bounded for some $ p\in (1,\infty)\setminus \{ 2 \} $ if and only if  $\Omega \in \mathcal{D}. $ As far as we know, the question of weather the boundedness of $\Pi_\Omega^+ $ on $ L^2(\Omega dA) $ implies that $ \Omega \in \mathcal{D} $ is still open.

In view of the above results in the unit disc, the following theorem comes as a surprise. 

\begin{thm}\label{thm:BergmanProjBounded}
	Let $ \omega $ a positive radial weight on a radial tree $X $. Then, for $ p\in (1,\infty) $ the Bergman projection $ P_\omega : L^p(\omega) \to L^p(\omega) $  is bounded and of weak type $(1,1)$.  
\end{thm}

It should also be noted that the source of this discrepancy between the theory in the unit disc and on trees is not completely understood. 
It is well understood that obstructions coming from local behaviour of kernels are absent in the discrete setting. Nonetheless, in this particular case the \eqref{cond:Dostanic} condition does not seem to be a local phenomenon, rather a global one. In fact its necessity is obtained by tasting the boundedness of the operator on functions of the form $ z^p \overline{z}^q $. 
It would therefore be desirable to get a better understanding of this phenomenon, that at the moment is still elusive. 

Regarding the positive Bergman projection $ P_\omega^+ $, in this case its definition is slightly less obvious. Our definition is in fact coherent with the definition used in the literature on trees \cite{DeMari2025, DeMari2024, CondeAlonso2026, CondeAlonso2024}. For the precise definition see Definition \ref{def:positiveBergman}. 

In the case of radial trees we have a complete characterization of weights which render the integral operator $ P_\omega^+ $ bounded on $ L^p(\omega) $. In order to do so we we need a few more notions from tree geometry.
Given $x\in X $ we denote by $T_x $ the tent at $x $, that is the set of vertices $y\in X $ such that 
the geodesic connecting $ o $ and $y $ passes from $x $. By $U_x $ we denote the pruned tent, i.e. $ U_x = T_x \setminus \{ x \}  $. 
We denote by $ s(v) $ the set of sons of $ v $, i.e. the set of vertices $u $  in $ T_v $ such that $ |u| = |v| + 1$. If $v=o$, then $s(o)$ is the set of all the neighbors. When $v\neq o$, there is a unique neighbor that do not lie in $s(v)$, we call it predecessor of $v$ and denote it by $p(v)$. We can see the predecessor as a function $p\colon X\setminus\{o\}\to X$ and in this way we iterate it to a function $p^k\colon X\setminus B(o,k-1)\to X$ by $p^k(v)=p(p^{k-1}(v))$.
 The family of all tents will be denoted by $ \mathcal{T} $ and the collection of all pruned tends by $ \mathcal{U} $. Furthermore if $ Q=T_x\in \mathcal{T} $ and $ n\in \mathbb{N}_0 $ we will denote by $ Q^{(n)} $ the tent with base point the $ n $-th predecessor of $x $. 
For $ U \in \mathcal{U} $ the set $ U^{(n)} $ is defined in a similar way.
We will also denote by $ B(x,N) $ the open ball with respect to the distance $d $ centred   at $x $ with radius $N $ . 

We will say that a weight $\omega  $ satisfies the lower doubling condition and we will write $ \omega \in \check{\mathcal{D}}_{\text{tree}} $  if there exists $\eta >0 $ and $ N>0 $ such that 
\begin{equation}
	\label{cond:Pplus} \eta \omega(T_x) \leq \omega (T_x \cap B(x,N)), \,\, \text{for all}\,\, x\in X.
\end{equation}

Here for $E \subseteq X $ we define $ \omega(E) = \sum_{x\in E} \omega(x)  $.
\begin{thm}\label{thm:PositiveBergman}
	Let $\omega $ be a positive radial weight on a radial tree $X $. Then, if $\omega$ satisfies condition $ \check{\mathcal{D}}_{\text{tree}} $   then the operator $ P_\omega^+ : L^p(\omega) \to L^p(\omega) $ is bounded for all $ p\in (1,\infty) $. Conversely, if $ P_\omega^+ : L^p(\omega) \to L^p(\omega) $ is bounded for some $ p\in (1,\infty) $ then $\omega$ satisfies condition $ \check{\mathcal{D}}_{\text{tree}}.$ 
\end{thm}

The techniques used to prove Theorem \ref{thm:BergmanProjBounded} and Theorem \ref{thm:PositiveBergman} are quite different. In the case of Theorem \ref{thm:BergmanProjBounded} we prove that the kernel of the operator satisfies the Calderon-Zygmund type conditions on trees introduced in \cite{DeMari2024}. Once the conditions on the kernel of $P_\omega $ are verified then the boundedness properties follow  from \cite[Theorem C]{DeMari2024}. 

On the other hand in order to establish the boundedness properties of the operator $P_\omega^+ $ for lower doubling weights $\omega $ we establish a sparse domination bound in the spirit of \cite{CondeAlonso2026, CondeAlonso2024}.  

The article is structured as follows. Section 2 is devoted to introducing of the main objects we need in the discrete setting and deriving an explicit formula for the reproducing kernel of the Bergman space. In Section 3, his formula will be useful in deriving kernel estimates and, consequently, boundedness results for the projector $P_\omega$. The study of the boundedness of $P_\omega^+$ involves different techniques and it is fully addressed in Section 4. Section 5 contains a discussion of the endpoint case $p=\infty$, where a suitable Bloch space is introduced. Finally in the Appendix we provide a short proof of the fact that the Dostani\'c condition implies $ \widehat{\mathcal{D}} $.  

Furthermore, throughout the rest of the paper we will use the following notation. If $ f,g  $ are non-negative functions on some set $ A $ we will write $ f \lesssim  g$ if there exists a constant $C>0 $ such that $ f(a) \leq C g(a) $ for all $ a \in A $. We will write $ f\asymp g $  if $ f\lesssim g $ and $ g \lesssim f $.  Finally, if $\omega$ is a positive measure on a space $X$, $Q$ is a subset of $X$ and $f\colon X\to \mathbb R$, we set 
\[\langle f\rangle_Q=\frac 1{\omega(Q)} \sum_{x\in Q }f(x)\omega(x),\]
for the average of $f$ on $Q$.

\section{Bergman projections on radial trees}

In this section we will derive an explicit expression for the kernel of the Bergman projection $P_\omega $ on a radial tree $ X $ . We also consider the finite dimensional space 
\[W_v = \bigg\{ \psi \in L^2(s(v)) : \sum_{u\in s(v)} \psi(u) = 0 \bigg\},\]
where the Lebesgue spaces are taken with respect to the counting measure. The dimension of $W_v $ is $ q(v)  $ and it is a reproducing kernel Hilbert space. 
Hence, the orthogonal projection $p_v: L^2(s(v)) \to W_v $ can be written as an integral operator with kernel given by the reproducing kernel of $ k_v $  of $W_v $. A straightforward verification shows that 
\begin{equation} \label{eq:LocalKernel}k_v(u,w) = \delta_{uw}-\frac{1}{q(v)}, \quad u,w \in s(v). \end{equation}
Indeed for $w\in s(v) $ we have that $k_v(\cdot,w) $ has zero mean and that, for $\psi \in W_v $,
\[ \langle \psi, k_w \rangle_{W_v} = \sum_{u\in s(v)}\psi(u)k_v(u,w)  = \psi(w)- \frac{1}{q(v)} \sum_{u\in s(v)} \psi(w) = \psi(w).  \]

For each $m,k \in \mathbb{N}_0, $ we consider the function 
\[ \varphi_k(m) = \sum_{j=0}^{m-k} \prod_{\ell=k}^{m-j-1} q(\ell)^{-1}, \]
with the convention that $ \varphi_k(m) = 0 $ if $ k > m $. It follows then by hypothesis \eqref{hyp:H0} that $ 1 \leq \varphi_k(m) < 2$ whenever $ k \leq m $.  
By abusing slightly the notation, we sometimes write $ \varphi_k(v) $ in place of $\varphi_k(|v|) $ for $ v\in X $.  

We can now define the extension operator $ E_v: W_v \to L^2(\omega) $
\[ E_v \psi (x) = \varphi_{|v|+1}(x)\sum _{u\in s(v)} \psi(u)  \mathbbm{1}_{T_u}(x), \quad x \in X. \]
This operator extends a given function $ \psi \in W_v $ to a function defined on the whole tree supported in $ U_v. $ 

To compute the adjoint of $ E_v $ let $ \psi \in W_v $ and $f \in L^2(\omega) $ we then have 
\begin{equation} \label{eq:ProjEv} 
	\begin{split}
		\langle E_v \psi, f  \rangle_{L^2(\omega)} &= \sum_{u \in s(v)} \psi(u) \sum_{x\in T_u} \varphi_{|v|+1}(x) f(x) \omega(x) \\
		&=\langle \psi , E_v^\times f \rangle_{L^2(s(v))}= \langle \psi , p_v E_v^\times f \rangle_{W_v}, 
	\end{split}
\end{equation}
where 
\[ E^\times_v f (u) = \sum _{x\in T_u} \varphi_{|v|+1}(x) f(x) \omega(x). \] 
Therefore $ E_v^* = p_v E^\times_v. $ 

\begin{prop}
	The image of the operator $E_v $  is a finite dimensional subspace of  $ A^2_\omega $. Furthermore, if $ u $ is any point in $s(v) $  we have that 
	\[ P_v = \Vert \varphi_{|v|+1} \mathbbm{1}_{T_u}\Vert^{-2} _{L^2(\omega)} E_v E_v^* \]
is the orthogonal projection of $L^2(\omega)  $ onto the subspace $ E_v(W_v) $.	
\end{prop}

\begin{proof}
	We start by proving that the operator $ P_v $ maps into $ A^2_\omega $. Since $ E_vf  $ vanishes on $X  \setminus U_v $, we have that it is harmonic in $ X \setminus T_v $. The harmonicity in $ U_v $ is guranteed by the fact that, for  $ m >0 $, 
	\begin{align*}
		q(m)&\varphi_k(m+1) + \varphi_k(m-1)  = q(m) \sum_{j=0}^{m+1-k} \prod _{\ell=k}^{m-j} q(\ell)^{-1} + \sum_{j=0}^{m-1-k} \prod_{\ell=k}^{m-j-2} q(\ell)^{-1} \\
		& = (q(m)+1) \sum_{j=0}^{m-1-k} \prod_{\ell=k}^{m-j-2} q (\ell)^{-1}  + q(m) \prod_{\ell = k}^{m-1}q(\ell)^{-1} + q(m) \prod_{\ell = k}^{m}q(\ell)^{-1} \\
		& = (q(m)+1) \sum_{j=0}^{m-1-k}  \prod_{\ell = k}^{m-j-2} q(\ell)^{-1} + (q(m)+1 ) \prod_{\ell=k}^{m-1}q(\ell)^{-1} \\
		& = (q(m)+1) \varphi_k(m).
	\end{align*}
Trivially the same holds for $ m = k $. 
	Finally we have that $E_v\psi $ is harmonic at $v $ by the fact that $ \varphi_k(k) =1  $ and that the mean of $\psi $ is zero over $s(v)$. Furthermore, since $\varphi_k $ is bounded and $\omega \in L^1(X)$, we also have that $E_v \psi \in A^2_\omega. $

	We now prove that $P_v $ is orthogonal projection $  L^2(\omega) $ onto $ E_v(W_v) $.  For $ u\in s(v) $ we have that 
	\[ E_v^\times E_v \psi (u) = \sum_{x\in T_u} \psi(u) \varphi_{|v|+1} (x)^2 \omega(x) = \Vert \varphi_{|v|+1 }\mathbbm{1}_{T_u}\Vert_{L^2(\omega)}^2 \psi(u), \]
	so that if $ f = E_v \psi \in E_v(W_v) $, then 
	\[ P_v f =  \Vert \varphi_{|v|+1 }\mathbbm{1}_{T_u}\Vert_{L^2(\omega)}^{-2} E_v p_v E^\times_v E_v \psi = E_v p_v \psi = f. \]
	Since $ P_v $ is self-adjoint with image contained is $E_v(W_v) $, we have the claim.
	 
\end{proof}

The importance of the family $\{ E_v(W_v) \}_{v \in X} $ of subspaces is provided by the following proposition.

\begin{prop}

	Let $X $ be a radial tree and let $\omega $ be a radial weight on $X $. Then we have the orthogonal decomposition 
	\[ A^2_\omega = \mathbb{C} \mathbbm{1}_X \oplus \bigoplus_{v\in X} E_v(W_v). \]
\end{prop}

\begin{proof}
	We start by observing that the subspaces are pairwise orthogonal. Clearly $\mathbbm{1}_X $ is orthogonal to $E_v(W_v) $ by the fact that for $ f\in A^2_\omega $ such that $  f(o)=0 $ and that $\omega $ is radial (we refer to~\cite{DeMari2024} and~\cite{CondeAlonso2024} for more details). Furthermore $ E_v(W_v)  $ and $E_w(W_w) $ are orthogonal by disjoint supports whenever $ v\wedge w \not\in \{ v,w \}. $ 
	Without loss of generality, it is sufficient to consider $v,w\in X $ such that $w\in U_v $. Let then $\psi \in W_v, \eta \in W_w $ . We have 

	\[ \langle E_v \psi, E_w \eta \rangle_{L^2(\omega)} = \psi(p^{|w|-|v|-1}(w)) \sum_{u\in s(w)} \eta(u) \sum_{x\in T_u} \varphi_{|v|+1}(x) \varphi_{|w|+1}(x) \omega(x) = 0, \]
	by the zero mean of $\eta $ and the radiality of the $\varphi_k$'s.
	
	Consider now $f\in A^2_\omega $, we define 
	\begin{align*} 
		\psi_{f,o}(u) & = f(u) - f(o) \\
		\psi_{f,v}(u) &  = f(u) - \frac{q(v)+1}{q(v)}f(v) + \frac{1}{q(v)}f(p(v)).
	\end{align*} 
By the fact that $f $ is harmonic, we have that $\psi_{f,v}\in W_v $ for each $v\in X $ . We put 
\[ g = f(o) + \sum_{v\in X} E_v \psi_{f,v} \]
and we prove that $f = g. $ The fact that $f(o)=g(o) $ comes again from the support of the functions $E_v(\psi) $ , while for $x\in X \setminus \{ o \} $ we have 
\begin{align}\label{eq:Identity}
	&g(x)  = f(o) + \sum_{\ell=1}^{|x|} E_{p^\ell(x)}\psi_{f,p^\ell(x)}(x) \notag \\
	     & = f(o)(1-\varphi_1(x)) + f(p^{|x|-1}(x)) \varphi_1(x) \notag \\
	     & \quad + \sum_{\ell=1}^{|x|-1} \Big[ f(p^{\ell-1}(x)) - \frac{q(p^\ell(x))+1}{q(p^\ell(x))}f(p^\ell(x)) + \frac{1}{q(p^\ell(x))}f(p^{\ell+1}(x)) \Big] \varphi_{|x|-\ell+1}(x) \notag \\
	     &=\sum_{\ell=1}^{|x|}a_\ell(x)f(p^\ell(x))+f(x),
\end{align}
with
\begin{align*}
	a_1(x)&=\varphi_{|x|-1}(x) - \frac{q(p(x))+1}{q(p(x))}, \\
	a_\ell(x)&=\varphi_{|x|-\ell}(x) - \frac{q(p^\ell(x))+1}{q(p^\ell(x))} \varphi_{|x|-\ell+1}(x) + \frac{1}{q(p^{\ell-1}(x))}\varphi_{|x|-\ell+2}(x),\\
	a_{|x|}(x)&=1 - \varphi_1(x) + \frac{1}{q(p^{|x|-1}(x))}\varphi_2(x).
\end{align*}
Now observe that by definition of $\varphi_k $ we have that 
\[ \varphi_k = q(k-1)(\varphi_{k-1}-1). \]
So we have that $a_\ell(x)=0$ for every $\ell=1,\dots,|x|$.
By putting together with  identity \eqref{eq:Identity}  we have that $f(x) = g(x) $. 
\end{proof}

The following corollary gives an explicit description of the kernel of the Bergman projection $P_\omega $. 

\begin{coro} \label{coro:decomp}
	We have that the orthogonal projection $ P _\omega $ can be written as 
\[ P_\omega f = \langle f  \rangle_X + \sum_{v\in X} P_vf. \]
\end{coro}

\noindent Notice that the operators $ P_v $ are integral operators, with kernels 
\begin{equation}\label{eq:KernelVertex} K_v(x,y) = \sum_{u\in s(v)} \sum_{w\in s(v)} k_v(u,w) \frac{\varphi_{|v|+1}(y) \varphi_{|v|+1}(x)}{\Vert \varphi_{|v|+1} \mathbbm{1}_{T_u}\Vert^2_{L^2(\omega)}} \mathbbm{1}_{T_w}(y) \mathbbm{1}_{T_u}(x),\qquad x,y \in X. 
\end{equation}
Therefore it follows by Corollary \ref{coro:decomp} that the kernel $ K_\omega $ of $ P_\omega $ is given by \begin{align*}\label{eq:BergmanKer}
	K_\omega(x,y) & = \frac{1}{\omega(X)} + \sum_{\ell=0}^{|x\wedge y|} K_{p^\ell(x \wedge y)}(x,y).
\end{align*}

Finally we conclude this section by defining the positive Bergman projection. 

\begin{defn} \label{def:positiveBergman}
	Let $ \omega \in L^1(X) $ a positive radial weight. For $f $ compactly supported we define the positive Bergman projection as the operator 
	\[ P_\omega^+ f (x) = \langle f \rangle_X + \sum_{v\in X} \sum_{y\in X} |K_v(x,y)| f(y) \omega(y). \]
	We also define the positive Bergman kernel as
	\[ K_\omega^+(x,y) = \frac{1}{\omega(X)} + \sum_{v\in X} |K_{v}(x,y)|. \]
\end{defn}

\section{Kernel estimates and the proof of Theorem \ref{thm:BergmanProjBounded}}

This section contains estimates for the kernels $k_v, K_v  $ and $ K_\omega $. After the desired estimates will have been obtained, Theorem \ref{thm:BergmanProjBounded} is going to be a direct application of the following theorem \cite[Theorem C]{DeMari2024}. 

 \begin{thm}\label{thm: C}

	Let $\omega $ be a finite positive weight on a tree $X $. Let $T $ be an integral operator bounded on $L^2(\omega) $ with a symmetric kernel $K $. Assume that $K $ satisfies the following conditions 	
	\begin{equation}
		\tag{H\"or}\label{cond:CalderonZ} 
		\sup _{Q \in \mathcal{T}} \sup_{x,y\in Q} \sum_{z\in X \setminus Q} |K_\omega(z,x) - K_\omega(z,y) | \omega(z) < + \infty 
	\end{equation}
	and 
	\begin{equation}
		\tag{S} \label{cond:Size} \sup_{Q\in \mathcal{T}}\sup_{x\in Q } \sum_{x\in Q^{(1)}\setminus Q}|K_\omega(z,x)| \omega(z) < + \infty.
	\end{equation}
	Then $ T $ is of weak type $(1,1) $ and hence of strong type $ (p,p) $ for every $ p \in (1,\infty) $. \end{thm}

The first lemma is an estimate of the size of the kernels. 

\begin{lem}
\label{prop: stime nuclei}The following estimates for the involved
kernels hold true. Let $v\in X$, we have
\begin{align}
|k_{v}(u,w)| & \asymp\begin{cases}
1 & u=w,\\
q(v)^{-1} & u\neq w,
\end{cases} &  & u,w\in s(v); \\
|K_{v}(x,y)| &  \asymp\begin{cases}
 \omega(U_{x\wedge y})^{-1}, & x\wedge y=v,\\
 \omega(T_{u})^{-1}, & x\wedge y\neq v,\,u\in s(v),
\end{cases} &  & x,y\in U_{v};\label{eq: stime K_v}\\
|K_{\omega}(x,y)| &  \lesssim\begin{cases}
 \omega(T_{x\wedge y})^{-1}, & x\wedge y\in\left\{ x,y\right\} ,\\
 \omega(U_{x\wedge y})^{-1}, & x\wedge y\not\in\left\{ x,y\right\} ,
\end{cases} &  & x,y\in X. 
\end{align}
\end{lem}

\begin{proof}
Let $v\in X$. The estimate for the kernel $k_{v}$ directly come
from the explicit formula \eqref{eq:LocalKernel}. Notice that, whenever
$x\wedge y=v$, we have that there exists $u, w\in s(v)$ with $ u\neq w $ 
such that $x\in T_{u}$ and $y\in T_{w}$. In this case, we have $k_{v}(u,w)=-q(v)^{-1}$
and hence by \eqref{eq:KernelVertex} 
\[
	|K_{v}(x,y)| = \frac{1}{q(v) \Vert \varphi_{|v|+1} \mathbbm{1}_{T_v}\Vert^2_{L^2(\omega)}} \asymp\frac{1}{q(v)\omega(T_{u})}=\frac{1}{\omega(U_{v})}.
\]
On the other hand, if $x\wedge y\neq v$, then $x,y\in T_{u}$ for
a certain $u\in s(v)$, so that $k_{v}(u,u)\asymp1$
and $K_{v}(x,y)\asymp\omega(T_{u})^{-1}$. Finally, we recall that
\[
	K_{\omega}(x,y)=\frac{1}{\omega(X)}+\sum^{|x\wedge y|}_{\ell=0}K_{p^{\ell}(x\wedge y)}(x,y).
\]
If $x\wedge y\in\left\{ x,y\right\} $, then the term for $\ell=0$
vanishes since $(x,y)\not\in {\rm supp}K_{x\wedge y}=U_{x\wedge y}\times U_{x\wedge y}$.
By using the previous estimates we get
\begin{align*}
 |K_{\omega}(x,y)| & \lesssim\begin{cases}
{\displaystyle \sum^{|x\wedge y|}_{\ell=1}|K_{p^{\ell}(x\wedge y)}(x,y)|,} & x\wedge y\in\left\{ x,y\right\} ,\\
{\displaystyle \sum^{|x\wedge y|}_{\ell=0}|K_{p^{\ell}(x\wedge y)}(x,y)|,} & x\wedge y\not\in\left\{ x,y\right\} 
\end{cases} \\
					     & \asymp 
\begin{cases}
{\displaystyle \sum^{|x\wedge y|}_{\ell=1}\frac{1}{\omega(T_{p^{\ell-1}(x\wedge y)})},} & x\wedge y\in\left\{ x,y\right\} ,\\
{\displaystyle \frac{1}{\omega(U_{x\wedge y})}+\sum^{|x\wedge y|}_{\ell=1}\frac{1}{\omega(T_{p^{\ell-1}(x\wedge y)})},} & x\wedge y\not\in\left\{ x,y\right\} .
\end{cases}
\end{align*}
We conclude by observing that, given a sector $T_x \in \mathcal{T} $, we have that 
\begin{align*}
	\omega(T_{p(x)}) = \omega(x) + q(x) \omega(T_x) \geq 2 \omega (T_x).
\end{align*}
Hence we have
\[
\frac{1}{\omega(U_{x\wedge y})}+\sum^{|x\wedge y|-1}_{\ell=0}\frac{1}{\omega(T_{p^{\ell}(x\wedge y)})}\leq\frac{1}{\omega(U_{x\wedge y})}+\frac{1}{\omega(T_{x\wedge y})}\sum^{|x\wedge y|-1}_{\ell=0}2^{-\ell}\lesssim\frac{1}{\omega(U_{x\wedge y})},
\]
and similarly in the other case.
\end{proof}

\begin{lem}\label{lem:CZ}
	Let $ \ell \in \mathbb{N}_0 $ and  $x,z,v \in X  $ such that, $x\in U_v $ and $ z\not\in T_v $. Then we have 
	\[ |K_{p^\ell(v\wedge z)}(z,x)-K_{p^\ell(v\wedge z)}(z,v)| \lesssim  \begin{cases} \displaystyle{ \frac{1}{\omega\big(T_{v\wedge z}^{(\ell-1)}\big)} \prod_{i=|v\wedge z| - \ell + 1 }^{|v|} q(i)^{-1}}, & \text{if }1\leq \ell \leq |v\wedge z|, \\
	\displaystyle{ \frac{1}{\omega(U_{v\wedge z})} \prod_{i=|v\wedge z| +1}^{|v|} q(i)^{-1} }, &  \text{if }\ell =0\text{, }z\neq z \wedge v.
	\end{cases} \]
\end{lem}
Notice that in the case $\ell=0$, we need to exclude the case $z=z\wedge v$, otherwise the left hand side of the inequality vanishes since $(z,v)\not\in {\rm supp}K_{p^{\ell(v\wedge z)}}={\rm supp}K_z$.
\begin{proof}  Let $x\in U_{v}$ and $z\not\in T_{v}$. Since $x\wedge z=v\wedge z=p^{n}(v)$
for some $n\in\left\{ 1,\dots,|v|\right\} $, we have that for $ \ell \geq 1 $ 
\begin{align*}
	&|K_{p^\ell(v \wedge z)}(z,x)-K_{p^\ell(v\wedge z)}(z,v)|  = |K_{p^{\ell+n}(v)}(z,x)-K_{p^{\ell+n}(v)}(z,v)| \\ 
							       & = k_{p^{\ell+n}(v)}(p^{\ell+n-1}(v),p^{\ell+n-1}(v)) \varphi_{|v|-n-\ell+1}(z)\frac{|\varphi_{|v|-n-\ell+1}(x)-\varphi_{|v|-n-\ell+1}(v)| }{\|\varphi_{|v|-n-\ell+1} \mathbbm{1}_{T_v^{(n+\ell-1)}}\|^{2}_{L^{2}(\omega)}}\\
 & \lesssim\frac{1}{\omega\big(T^{(n+\ell-1)}_{v}\big)}|\varphi_{|v|-n-\ell+1}(x)-\varphi_{|v|-n-\ell+1}(v)|.
\end{align*}
Furthermore for $ m\geq k $ we have that  
\begin{align*}
\varphi_{k}(m+n)-\varphi_{k}(m) & =\sum^{m+n-k}_{j=0}\prod^{m+n-j-1}_{i=k}q(i)^{-1}-\sum^{m-k}_{j=0}\prod^{m-j-1}_{i=k}q(i)^{-1}\\
&=\sum^{n-1}_{j=0}\prod^{m+n-j-1}_{i=k}q(i)^{-1}\\
 & =\prod^{m}_{i=k}q(i)^{-1}\sum^{n-1}_{j=0}\prod^{m+n-j-1}_{i=m+1}q(i)^{-1}\lesssim\prod^{m}_{i=k}q(i)^{-1}.
\end{align*}
When $\ell =0 $ and $ z\neq z\wedge v $  there are $u,w\in s(z\wedge v)$
such that $z\in T_{u}$ and $v\in T_{w}$. Hence, in this case
\begin{align*}
	 |K_{v\wedge z}(z,x) - K_{v\wedge z}(z,v)| &=	|k_{p^{n}(v)}(u,w)|\varphi_{|v|-n+1}(z)\frac{|\varphi_{|v|-n+1}(x)-\varphi_{|v|-n+1}(v)| }{\|\varphi_{|v|-n+1} \mathbbm{1}_{T_u}\|^{2}_{L^{2}(\omega)}} \\ 
& \lesssim \frac{1}{\omega\big(U^{(n)}_{v}\big)}\prod^{|v|}_{i=|v|-n+1}q(i)^{-1},
\end{align*}
by using that $|k_{p^{n}(v)}(u,w)|=q(p^{n}(v))^{-1}$.
\end{proof}

Theorem \ref{thm:BergmanProjBounded} will be an immediate consequence of the following proposition. 

\begin{prop} \label{prop:HoS}
	Let $ \omega $ a positive radial weight on a radial tree $X $. Then the kernel $ K_\omega $ satisfies \eqref{cond:CalderonZ} and \eqref{cond:Size}.  
\end{prop}

\begin{proof}
	We start by noticing that, by triangular inequality, it is sufficient to prove condition \eqref{cond:CalderonZ} with $ Q = T_v$, $x\neq v $, and $ y =v. $ Let $x\in U_v $ and $z\not\in T_v $. Let also $n\in \mathbb{N} $ such that $x\wedge z = v\wedge z = p^n(v) $. We then have  
	\begin{align*}
		|K_\omega(z,x)-K_\omega(z,v)| & \leq |K_{v\wedge z}(z,x) - K_{v\wedge z}(z,v)| + \sum_{\ell=1}^{|v|-n} |K_{p^\ell(v\wedge z)}(z,x) - K_{p^\ell(v\wedge z)}(z,v)| \\
				& =\text{I}+\text{II} .
	\end{align*}
	The term $\text{I}$ vanishes for $z=z\wedge v$ and is controlled by Lemma~\ref{lem:CZ} otherwise. Notice that by Lemma \ref{lem:CZ} it follows that 
	\begin{align*}
		\,\, \text{II}\,\,& \lesssim \sum_{\ell=1}^{|v|-n} \frac{1}{\omega\big(T_v^{(n+\ell-1)}\big)}\prod _{i=|v|-n-\ell+1}^{|v|}\frac{1}{q(i)} \\
				  & = \frac{1}{\omega\big(T_v^{(n)}\big)} \prod_{i=|v|-n}^{|v|} \frac{1}{q(i)} \sum_{\ell=1}^{|v|-n} \frac{\omega\big(T_v^{(n)}\big)}{\omega\big(T_v^{(n+\ell-1)}\big)} \prod_{t=1}^{\ell-1}\frac{1}{q(|v|-n-\ell+t) } \\
				  & \lesssim  \frac{1}{\omega\big(T_v^{(n)}\big)} \prod_{i=|v|-n}^{|v|}\frac{1}{q(i)}\lesssim \frac 1{\omega\big(U_v^{(n)}\big)}\prod_{i=|v|-n+1}^{|v|}\frac 1{q(i)},
	\end{align*}
	where we controlled with the same estimate used for $\text{I}$.
By putting everything together, we have
\begin{align*}
 & \sum_{z\in X\setminus T_{v}}\left|K_\omega(z,x)-K_\omega(z,v)\right|\omega(z)\\
= & \sum^{\left|v\right|}_{n=1}\left[\left|K_\omega(p^{n}(v),x)-K_\omega(p^{n}(v),v)\right|\omega(p^{n}(v))+\sum_{z\in U^{(n)}_{v}\setminus T^{(n-1)}_{v}}\left|K_\omega(z,x)-K_\omega(z,v)\right|\omega(z)\right]\\
\lesssim &\sum_{n=1}^{|v|}\Bigg[\frac{\omega(p^n(v))}{\omega\big(T_v^{(n)}\big)}\prod_{i=|v|-n}^{|v|}\frac 1{q(i)}+\frac{\omega\big(U_v^{(n)}\setminus T_v^{(n-1)}\big)}{\omega\big(U_v^{(n)}\big)}\prod_{i=|v|-n+1}^{|v|}\frac 1{q(i)}\Bigg]\\
\lesssim& \sum_{n=1}^{|v|}\prod_{i=|v|-n+1}^{|v|}\frac 1{q(i)}\lesssim 1.
\end{align*}
Similarly, we prove the size condition (\ref{cond:Size}) for $Q=T_{v}$
and $x=v$, since the other cases follows by applying the Hörmander
condition (\ref{cond:CalderonZ}) to the kernel. From Proposition \ref{prop: stime nuclei},
we have
\begin{align*}
\sum_{z\in T^{(1)}_{v}\setminus T_{v}}\left|K_\omega(z,v)\right|\omega(z) & =\left|K(p(v),v)\right|\omega(p(v))+\sum_{z\in U^{(1)}_{v}\setminus T_{v}}\left|K_\omega(z,v)\right|\omega(z) \\
								   & \lesssim\frac{\omega(p(v))}{\omega\big(T^{(1)}_{v}\big)}+\frac{\omega\big(U^{(1)}_{v}\setminus T_{v}\big)}{\omega\big(U^{(1)}_{v}\big)}\lesssim1.
\end{align*}
\end{proof}

\section{Sparse domination and the boundedness of $P_\omega^+ $ }

The proof of boundedness of the operator $P_\omega^+ $ goes through the sparse domination of the operator. Let $\omega $ be a positive weight on the tree and $\Lambda > 1 $. A family $\mathcal{S} $  of subsets will be called $ \Lambda $-Carleson in $(X,\omega) $ if for every $\Omega \subseteq X $, we have 
\[ \sum_{Q\in \mathcal{S}, Q\subseteq \Omega} \omega(Q) \leq \Lambda \omega(\Omega). \]
If the family $\mathcal{S} $ is nested, that is $ Q\cap R \neq \emptyset  $ for some $ Q,R \in \mathcal{S} $ implies either $ Q \subseteq R $ or $ R \subseteq Q $, it suffices to check the condition on $\Omega \in \mathcal{S} $. 

In the context of metric measure spaces without atoms, the notion of $\Lambda- $ Carleson is equivalent to the notion of sparseness~\cite{Verbitsky}. Here we present a modification of the argument in the case $\mathcal{S} = \mathcal{T}. $ 
The main observation, a simple one nonetheless, is that a weight $\omega$ renders the families of sets $ \mathcal{T} $ and $ \mathcal{U} $ $\Lambda$-Carleson for some $ \Lambda >1 $ if and only if $\omega $ is $ \check{\mathcal{D}}_{\text{tree}} $.

\begin{prop}
\label{prop: equiv carl}Given $\omega\in L^{1}(X)$ be positive.
\begin{enumerate}
	\item[(a)] If $ \omega \in \check{\mathcal{D}}_{\text{tree}} $ then $ \mathcal{T} $ is $\Lambda $-Carleson family of sets in $ (X,\omega) $.
	\item[(b)] If $\mathcal{T} $ is $\Lambda $-Carleson in $(X,\omega) $ for some $\Lambda>1 $ , then for every $\eta\in(0,1)$
	there exists $N\in\mathbb{N}$ such that
	\begin{equation}
		\label{eq: equiv carl}
		\eta \omega(T_x) \leq \omega (T_x\cap B(x,N)), \quad x \in X.
	\end{equation}
	Hence $\omega\in\check{\mathcal D}_{\text{tree}}$.
\end{enumerate}
\end{prop}

\begin{proof}
	The condition $ \check{\mathcal{D}}_{\text{tree}} $  is clearly sufficient to have
that $ \mathcal{T} $ is Carleson in $(X,\omega) $   since for every $Q\in\mathcal{T}$
we have
\[
\sum_{x\in Q}\omega(T_{x})\leq\frac{1}{\eta}\sum_{x\in Q}\omega(T_{x}\cap B(x,N))=\frac{1}{\eta}\sum_{x\in Q}\sum_{y\in T_{x}\cap B(x,N)}\omega(y)\leq\frac{N}{\eta}\omega(Q).
\]
On the other hand, assume $\mathcal{T} $ is $\Lambda $-Carleson. Let
$x\in X$ and $\eta\in(0,1)$. We
set $N=\left\lfloor \frac{\Lambda}{1-\eta}\right\rfloor $
\begin{align*}
	\Lambda\omega(T_{x}) & \geq\sum_{y\in T_{x}}\omega(T_{y})=\sum_{y\in T_{x}}\sum_{z\in T_{y}}\omega(z)\\
 & =\sum_{z\in T_{x}}(d(x,z)+1)\omega(z)\geq\frac{\Lambda}{1-\eta}\omega(T_{x}\setminus B(x,N)).
\end{align*}
Hence we have
\[
\omega(T_{x}\cap B(x,N))=\omega(T_{x})-\omega(T_{x}\setminus B(x,N))\geq\eta\omega(T_{x}).
\]
\end{proof}

Notice that $\mathcal{T} $ is Carleson in $ (X,\omega) $ if and only if $\mathcal{U} $ is.  
Indeed, assume that $\mathcal{T} $ is Carleson in $ (X,\omega) $, given $x\in X$,
we have
\[
\sum_{\substack{S\in\mathcal{U}\\
S\subseteq U_{x}
}
}\omega(S)=\sum_{\substack{Q\in\mathcal{T}\\
Q\subseteq U_{x}
}
}\omega(Q)=\sum_{\substack{P\in\mathcal{T}\\
P^{(1)}=T_{x}
}
}\sum_{\substack{Q\in\mathcal{T}\\
Q\subseteq P
}
}\omega(Q)\leq\Lambda\sum_{\substack{P\in\mathcal{T}\\
P^{(1)}=T_{x}
}
}\omega(P)=\Lambda\omega(U_{x}).
\]
On the other hand, if $\omega\in\mathcal{C}_{\mathcal{U}}$ and $x\in X$,
then
\[
\sum_{\substack{S\in\mathcal{T}\\
S\subseteq T_{x}
}
}\omega(S)=\omega(T_x)+\sum_{\substack{S\in\mathcal{T}\\
S\subsetneq T_x
}
}\omega(S)=\omega(T_x)+\sum_{\substack{P\in\mathcal{U}\\
P\subseteq U_{x}
}
}\omega(P)\leq\omega(T_{x})+\Lambda\omega(U_{x})\leq(\Lambda+1)\omega(T_{x}).
\]
Hence, it follows that the previous characterization can be repeated
to characterize measures $\omega\in \check{\mathcal{D}}_{\text{tree}}$ by
replacing (\ref{eq: equiv carl}) with 
\begin{equation}
\eta\omega(U_{x})\leq\omega(U_{x}\cap B(x,N-1)),\qquad x\in X.\label{eq: equiv Carl U}
\end{equation}

We have the following characterization of the measure that make $P^{+}_{\omega}$
bounded on $L^{p}(\omega)$. The same argument can be found in the recent preprint~\cite{CondeAlonso2026}. Nevertheless, we include a direct proof, as the argument is particularly simple in our setting and highlights the role of the geometry of the tree and the decay of the measure.
\begin{thm}
Given $\omega\in L^{1}(X)$ positive. The following conditions are
equivalent.
\begin{enumerate}
	\item[(a)] $ \omega \in \check{\mathcal{D}}_{\text{tree}} $;
	\item[(b)] $P^{+}_{\omega}$ is bounded on $L^{p}(\omega)$ for
	some (hence for all) $1<p<\infty$.
\end{enumerate}
\end{thm}

\begin{proof}
Assume (i) holds. We prove a sparse domination result for $P^{+}_{\omega}$.
Indeed, it directly follows from the estimates on the kernel, for
$x\in X$ and $f\in L^p(\omega)$, $f(o)=0 $ we have
\begin{align*}
P^{+}_{\omega}f(x) & =\sum^{|x|}_{n=1}\left[K^{+}_{\omega}(x,p^{n}(x))f(p^{n}(x))\omega(p^{n}(x))
		 +  \sum_{y\in U^{(n)}_{x}\setminus T^{(n-1)}_{x}}K^{+}_{\omega}(x,y)f(y)\omega(y)\right] \\
		   & \quad  + \sum_{y\in T_{x}}K^{+}_{\omega}(x,y)f(y)\omega(y)\\
 & \leq\sum^{|x|}_{n=0}\frac{1}{\omega(T^{(n)}_{x})}\sum_{y\in T^{(n)}_{x}}|f(y)|\omega(y)+\sum^{|x|}_{n=1}\frac{1}{\omega(U^{(n)}_{x})}\sum_{y\in U^{(n)}_{x}}|f(y)|\omega(y) \\
 & =:\mathcal{A}_{\mathcal{T}}f(x)+\mathcal{A}_{\mathcal{U}}f(x).
\end{align*}
 We can then argue
as usual to show that $\mathcal{A}_{\mathcal{S}}$ is
bounded on $L^{p}(\omega)$ and of weak-type (1,1) estimate for $\mathcal{S} = \mathcal{T}, \mathcal{U} $. For completeness, we recall the proof of weak-type (1,1) for $\mathcal{A}_{\mathcal{T}}$ (the proof for $\mathcal A_{\mathcal U}$ is analogous). Take $f\in L^1(\omega)$. By a classical argument~\cite[Exercise 1.4.14]{Grafakos} it is sufficient to prove that 
\[ \sup_{\substack{F\subseteq X\\ 0<\omega(F)<\infty}}\inf_{\substack{F'\subseteq F\\ \omega(F')\geq \frac 12\omega(F)}}\sum_{x\in F'}\mathcal A_\mathcal Tf(x)\omega(x)\lesssim \|f\|_1,\]
as the LHS is equivalent to $\|\mathcal A_{\mathcal T}f\|_{1,\infty}$. Take $F\subseteq X$ with $0<\omega(F)<\infty$. Take $\eta\in(0,1)$. By Proposition~\ref{prop: equiv carl}, there exists $N\in \mathbb N$ depending on $\omega$ and $\eta$ such that~\eqref{eq: equiv carl} (and~\eqref{eq: equiv Carl U}) holds. We set
\[ F'=\bigg\{x\in F\colon \mathcal A_{\mathcal T}1_F(x)\leq \frac {8N}\eta\bigg\}.\]
For every $v\in X$, we define $E_{T_v}=T_v\cap B(v,N)$ (and $E_{U_v}=U_v\cap B(v,N-1)$ for the case of $\mathcal A_{\mathcal U}$). We first observe that
\begin{align*}
	\omega(F\setminus F')&<\frac \eta {8N}\sum_{x\in F}\mathcal A_{\mathcal T}1_F(x)\omega(x)= \frac \eta {8N}\sum_{Q\in \mathcal T}\langle 1_F\rangle_Q^2\omega(Q)\\
	&\leq \frac 1{8N}\sum_{Q\in \mathcal T}\sum_{x\in E_Q}|(\mathcal M_\mathcal T 1_F)(x)|^2\omega(x)\leq \frac 18 \|\mathcal M_{\mathcal T}1_F\|_2^2\leq \frac 12\omega(F),
\end{align*}
by using that the maximal function $\mathcal{M}_{\mathcal T}$ associated to sectors satisfies  $\|\mathcal{M}_{\mathcal T}f\|_2\leq 2\|f\|_2$ (see Doob's maximal inequality~\cite[§14.6]{Williams}). This implies that $\omega(F')\geq \omega(F)/2$. Let $x\in X$. If $Q_0\in\mathcal T$ is the smallest sector containing $x$ and such that $\omega(Q_0\cap F')>0$, then for every $y\in Q_0\cap F'$, we have that 
\[\mathcal A_\mathcal T 1_{F'}(x)\leq \mathcal A_\mathcal T 1_{F'}(y)\leq\frac {8N}\eta.\]
We conclude that
\[\sum_{x\in F'}\mathcal A_{\mathcal T}f(x)\omega(x)=\sum_{Q\in\mathcal T}\langle |f|\rangle_Q \omega(Q\cap F')= \sum_{x\in X} |f(x)|\mathcal A_{\mathcal T}1_{F'}(x)\omega(x)\lesssim \|f\|_1.\]
We have proved that $\mathcal A_{\mathcal T}\colon L^1(\omega)\to L^{1,\infty}(\omega)$ is bounded and the same holds for $P_\omega^+$. 
%
On the other hand, assume (ii) holds
for $1<p<\infty$. Let $Q\in\mathcal{T}$ and $f=g=\mathbbm{1}_{Q}$. Since
$f,g\in L^{r}(\omega)$ for every $1<r<\infty$ , we have
that 
\[
\left\langle P^{+}_{\omega}f,g\right\rangle \lesssim\left\Vert f\right\Vert _{p}\left\Vert g\right\Vert _{p'}=\omega(Q).
\]
Notice that 
\begin{align*}
\left\langle P^{+}_{\omega}f,g\right\rangle  & =\sum_{x\in Q}\sum_{y\in Q}K^{+}_{\omega}(x,y)\omega(y)\omega(x)\\
 & \geq\sum_{x\in Q}\bigg[\sum_{y\in T_{x}}K^{+}_{\omega}(x,y)\omega(y)+\sum^{|x|-|Q|}_{k=1}\sum_{y\in U^{(k)}_{x}\setminus T^{(k-1)}_{x}}K^{+}_{\omega}(x,y)\omega(y)\bigg]\omega(x)\\
 & \gtrsim\sum_{x\in Q}\bigg[1+\sum^{|x|-|Q|}_{k=1}\frac{\omega\big(U^{(k)}_{x}\setminus T^{(k-1)}_{x}\big)}{\omega\big(U^{(k)}_{x}\big)}\bigg]\omega(x)\\
 & =\sum_{x\in Q}\bigg[1+\sum^{|x|-|Q|}_{k=1}\Big(1-\frac{1}{q(p^{k}(x))}\Big)\bigg]\omega(x) \\
 & \gtrsim\sum_{x\in Q}\#\left\{ R\in\mathcal{T}\colon x\in R\subseteq Q\right\} \omega(x)=\sum_{\substack{R\in\mathcal{T}\\
R\subseteq Q
}
}\omega(R).
\end{align*}
By putting all together, we have that $\mathcal{T} $ is Carleson.
\end{proof}

\section{Endpoint estimates for the Bergman projection}

In this section we prove  $P_\omega$ the end point estimate $ P_\omega : L^\infty(X) \to \mathcal{B} $. In order to define a suitable Bloch space in the discrete context, we rely on the classic notion of discrete gradient $\nabla g(x)=g(p(x))-g(x)$ for $ x \in U_o $ . Hence we define the Bloch space as
\[ \mathcal B=\{f\colon X\to \mathbb R\colon \|f\|_{\mathcal B}:=|f(o)|+\|\nabla f \mathbbm{1}_{U_o}\|_{\infty}<+\infty \}.\]

As a direct corollary of Theorem \cite[Theorem C]{DeMari2025}, we have that $P_\omega$ is bounded from $L^\infty(X)$ to $BMO(\omega)$ and maps $ H^1(\omega) $ into $L^1(\omega)$, where the spaces $BMO(\omega)$ and $H^1(\omega)$ are defined as follows.
\[
\|f\|_{\rm BMO(\omega)} = \langle | f| \rangle_X + \sup_{x\in X}  |f(x)-\langle f\rangle_{T_x}|+ \sup_{\substack{Q\in \mathcal{T}\\ Q \neq X }} \langle |f-\langle f\rangle_{Q^{(1)}}|\rangle_{Q}.\]
 	In what follows we shall also use the norm equivalence (see~\cite{CondeAlonso2024}) given by
 	\begin{equation}\label{eq:equivnorm}
 		 		\|f\|_{\rm BMO(\omega)}\asymp \langle | f| \rangle_X+\sup_{Q\in\mathcal T}\langle  |f-\langle{f}\rangle_{Q}|\rangle_Q+\sup_{\substack{Q\in\mathcal T\\Q\ne X}} |\langle f\rangle_Q-\langle f\rangle_{Q^{(1)}}|.
 		 	\end{equation}
 	 	Finally, the dual space of ${\rm BMO(\omega)}$ is the atomic Hardy space defined by square function
 	 	\[H^1(\omega)=\bigg\{f\in L^1(\omega)\colon \|f\|_{H^1} =
 	 	\sum_{x\in X} \bigg(|f(x)-\langle f\rangle_{T_x}|^2+\sum_{\ell=1}^{|x|}|\langle f\rangle_{T_x^{(\ell-1)}}-\langle f\rangle_{T_x^{(\ell)}}|^2\bigg)^{\frac 12}\bigg\}.
 	 	\]
It is easy to see that ${\rm BMO(\omega)}\subseteq \mathcal B$. Indeed, for $f\in {\rm BMO(\omega)}$, we write
\[|f(x)-f(p(x))|\leq |f(x)-\langle f\rangle_{T_x}|+|\langle f\rangle_{T_x}-\langle f\rangle_{T_x^{(1)}}|+|\langle f\rangle_{T_x^{(1)}}-f(p(x))|\lesssim \|f\|_{\rm BMO(\omega)},\]
by applying the definition of ${\rm BMO(\omega)}$ to the first and the last term, and by using~\eqref{eq:equivnorm} for the middle one.

It turns out, that if $\omega$ is sufficiently decaying the two spaces are the same with equivalent norms. Nevertheless, the equality does not hold in general.

\begin{prop}
	Let $\omega\in L^1(X)$ be positive. Then the following are equivalent
	\begin{enumerate}
		\item[(a)] $\mathcal B={\rm BMO(\omega)}$ with equivalent norms.
		\item[(b)] $\omega \in \check{\mathcal D}_{\rm tree}$.
	\end{enumerate}
\end{prop}
\begin{proof}
	We first assume $\omega\in\check{\mathcal D}_{\rm tree}$. It is sufficient to prove $\mathcal B\subseteq {\rm BMO(\omega)}$. Take $f\in \mathcal B$. We notice that for $y\in T_x$ with $d(x,y)=n$, we can iterate the definition of $\|f\|_{\mathcal B}$ to obtain
	\[|f(x)-f(y)|\leq \sum_{j=1}^{n}|f(p^j(y))-f(p^{j-1}(y))|\leq n\|f\|_\mathcal B. \]
	Hence we can estimate
	\begin{align*}
		|f(x)-\langle f\rangle_{T_x}|&\leq \frac 1{\omega(T_x)}\sum_{y\in T_x}|f(x)-f(y)|\omega(y)\leq \frac {\|f\|_{\mathcal B}}{\omega(T_x)}\sum_{y\in T_x}d(y,x)\omega(y)\\
		&\leq \frac {\|f\|_{\mathcal B}}{\omega(T_x)}\sum_{R\in\mathcal \mathcal{T},R\subseteq T_x} \omega(R)\leq \Lambda\|f\|_{\mathcal B},
	\end{align*}
	by using the fact that $\mathcal T$ is $\Lambda$-Carleson in $(X,\omega)$.  Also, by similar computation,
	\begin{align*}
		\langle|f-\langle f\rangle_{Q^{1}}|\rangle_Q&\leq\frac 1{\omega(Q)}\frac 1{\omega(Q^{(1)})}\sum_{x\in Q}\sum_{y\in Q^{(1)}}|f(x)-f(y)|\omega(y)\omega(x)\\
		&\leq\|f\|_{\mathcal B}\frac 1{\omega(Q)}\frac 1{\omega(Q^{(1)})}\sum_{y\in Q^{(1)}}\sum_{x\in Q}d(x,y)\omega(x)\omega(y)\\
		&\leq\frac 1{\omega(Q)}\frac 1{\omega(Q^{(1)})}\sum_{R\in\mathcal \mathcal{T},R\subseteq Q^{(1)}}\omega(R)\sum_{y\in Q^{(1)}}\omega(y)\leq \Lambda \|f\|_{\mathcal B}.
	\end{align*}
	Finally, always by the same approach one also gets 
	\[\bigg|\sum_{x\in X}f(x)\omega(x)\bigg|\leq|f(o)|(\omega(X)-\omega(o))+\sum_{y\neq o}|f(y)-f(o)|\omega(y)\lesssim \|f\|_{\mathcal B}. \]
	Hence $f\in {\rm BMO(\omega)}$. 
	It remains to prove that if $\omega$ is such that $\mathcal B={\rm BMO(\omega)}$, then there exists $\Lambda>1$ such that $\mathcal T$ is $\Lambda$-Carleson in $(X,\omega)$. Let $f(x)=|x|$. It is straightforward that $\|f\|_{\mathcal B}=1$ and so, by hypothesis $f\in {\rm BMO(\omega)}$. Hence by definition of ${\rm BMO(\omega)}$-norm we have
	\begin{align*}
		\|f\|_{\rm BMO(\omega)}&\geq|f(x)-\langle f\rangle_{T_x}|\\
		&=\bigg||x|-\frac 1{\omega(T_x)}\sum_{y\in T_x}(|x|+d(x,y))\omega(y)\bigg|=\sum_{R\in\mathcal T,R\subseteq T_x}\omega(R),
	\end{align*}
	which proves that $\mathcal T$ are $\Lambda $-Carleson in $(X,\omega)$ with $ \Lambda = \|f\|_{\rm BMO(\omega)} $ .
\end{proof}

\begin{coro} 
	Let $\omega $ a positive radial weight on $X $. Then $ P_\omega : L^\infty(\omega) \to \mathcal{B} $ is bounded.
\end{coro}

\section{Appendix} \label{sec: Appendix}

Let $\Omega \in L^1(0,1) $ be a weight function. We will assume that $\Omega $ does not vanish identically on any interval $(r,1)$, $0<r<1 $.  
We will denote by $ \widehat{\Omega}(r): = \int_r^1 \Omega(t)dt $ the cumulative weight function. 
We will also use the notation 
 \[ M_\Omega(x) = \int_0^1 \Omega(t) t^{x+1} dx. \]

 We would like to prove that if a weight satisfies the Dostani\'c condition for some $p\neq 2 $ then it satisfies the upper doubling condition. We will assume without loss of generality that $p\in (1,2) $.

For a function $f $  of one variable we will use the following notation 
\[  \Delta(f;x,y) = \frac{f(y)-f(x)}{y-x} \]
whenever the increment quotient is well defined. 
 We start with the following elementary lemma about concave functions. 

\begin{lem}
	Let $f:(0,\infty)\to \mathbb{R} $ a non-negative non-decreasing concave function such that $f(x) = o(x) $ as $x \to \infty. $ Suppose, furthermore, that there exists $\alpha\in(\frac{1}{2},1) $ and $C>0 $  such that for all $x>0 $ we have
	\[ \Delta(f;x,\alpha x) \leq \Delta\Big(f;\frac{\alpha}{2\alpha-1}x,x\Big) + \frac{C}{x}.   \]
Then for all $x>0 $ we have
\[ f(x)-f(\alpha x) \leq C. \]
\end{lem}

\begin{proof}
	Let $x>0 $, we have that
\begin{align} \label{eq:recursion}
		 \begin{split}
		 	\Delta(f;x,\alpha x)  & \leq \Delta\Big(f; \frac{\alpha}{2\alpha-1}x,x\Big) + \frac{C}{x}
		 	\\&= \Delta\Big(f; \frac{\alpha^2}{2\alpha -1} \frac{x}{\alpha}, x \Big) + \frac{C}{x} 
		 	\leq \Delta\Big(f; \frac{x}{\alpha},x\Big) + \frac{C}{x}. 
		 \end{split}
\end{align}
The only non trivial step in the above chain of inequalities is the last one where we have used the concavity of $f $ and the fact that $ \frac{\alpha^2}{2\alpha -1 }>1 $. 
By iterating $n $ times the same inequality we that 
\begin{equation}\label{eq:complete_rec} \Delta(f;x,\alpha x) \leq \Delta\Big(f; \frac{x}{\alpha^n}, \frac{x}{\alpha^{n-1}}\Big) + \frac{C}{x} \sum_{i=0}^{n-1} \alpha^i.  
\end{equation}
Notice that 
\[ 0\leq \Delta\Big(f; \frac{x}{\alpha^n},\frac{x}{\alpha ^{n-1}}\Big) =  \frac{f(\frac{x}{\alpha ^n})-f(\frac{x}{\alpha^{n-1}})}{\frac{x}{\alpha ^{n}}(1-\alpha)} \leq \frac{1}{1-\alpha} \frac{\alpha ^n f(\frac{x}{\alpha^n})}{x} = o(1) \]
as $n\to \infty $. 
 Therefore taking the limit as $n\to \infty $ in inequality \eqref{eq:complete_rec} we find that 
\[ \Delta(f;x,\alpha x) \leq \frac{C}{x(1-\alpha)} \]
which is equivalent to the desired inequality.
\end{proof}

We now show how to apply the previous lemma to the function $f(x) = -\ln M_\Omega(x) $.
Without loss of generality we can assume that $\widehat{\Omega}(0) = 1 $ so that $f $ is non-negative. The concavity of $f $ follows from H\"older's inequality. The fact that $f(x) = o(x) $ as $x\to \infty $ follows from the fact that $\Omega $ does not vanish identically on any interval $(r,1), r\in (0,1) $. Indeed,  
\[ \frac{f(x)}{x} = - \ln\bigg(\Big( \int_0^1 t^x \Omega(t) t dt  \Big)^{\frac{1}{x}}\bigg). \]
The quantity inside the logarithm, as $x \to \infty $,  tends to the essential supremum of $t $ with respect to the measure $\Omega(t) t dt $. Since $\Omega $ does not vanish a.e. on any interval $(r,1) $ we have that the essential supremum is equal to $1 $. 
Finally the Dostani\'c condition for $p\in (1,2) $ is equivalent to the inequality
\begin{equation}
	f(2n) - \frac{1}{p} f(pn) - \frac{1}{p'} f(p'n) \leq C, \qquad n\in \mathbb{N}.
\end{equation}
Notice that by the concavity of $f $ we have that for any $x >0 $ 
\[ f(x+1) - f(x) = \Delta(f;x,x+1) \leq \Delta(f;0,1). \]
The last quantity is a constant independent of $x $, therefore the Dostani\'c condition implies that 
\[ f(x)-\frac{1}{p}f\Big(\frac{p}{2}x\Big)-\frac{1}{p'}f\Big(\frac{p'}{2}x\Big) \leq C', \qquad x>0. \]
Rearranging the terms we find that
\[ \frac{2-p}{2p}x \Big(\Delta\Big(f;x,\frac{p}{2}x\Big)-\Delta\Big(f;x,\frac{p'}{2}x\Big)\Big) \leq C',\qquad x>0 \]
	which is equivalent to the hypothesis of the previous lemma with $\alpha = \frac{p}{2} $ and $C = \frac{p}{2-p}C' $. Therefore we can conclude that
\[ f(x)-f\Big(\frac{p}{2}x\Big) \leq C, \qquad  x>0. \]
This is equivalent to the upper doubling condition \cite[Lemma 7]{Pelaez2024}.

\section{Use of Artificial Intelligence}

No artificial intelligence tool has been used to develop the theorems or to search the literature or write the main part of the article except for the appendix. 
The theorem in the appendix was first proved by Chat GPT 5.5 in an involved way. The proof presented in the appendix is an elaboration and simplification, by the authors, without the use of AI of the Chat GPT proof. Finally typos and minor errors were corrected with the help of Claude AI.

\nocite{*}
\bibliographystyle{plain}
\bibliography{biblio}

\end{document}